\documentclass[11pt,a4paper,reqno]{amsart}
\usepackage{latexsym}
\usepackage{amssymb,amsthm,amsmath}
\usepackage[dvips]{graphicx}
\usepackage{xcolor}
\usepackage{mathtools} 
\usepackage{enumitem} 

\calclayout

\theoremstyle{plain}
\newtheorem{theorem}{Theorem}
\newtheorem{lemma}[theorem]{Lemma}
\newtheorem{corollary}[theorem]{Corollary}

\theoremstyle{definition}
\newtheorem{definition}[theorem]{Definition}
\theoremstyle{remark}
\newtheorem{remark}[theorem]{Remark}

\makeatletter	
\@namedef{subjclassname@2020}{%
  \textup{2020} Mathematics Subject Classification}
\makeatother

\def\endmark{\hskip 2em\begin{picture}(8,10)
\put(0,0){$\Box$} \put(2,0){\rule{1.9mm}{0.3mm}}
\put(6.5,0){\rule{0.3mm}{1.9mm}}
\end{picture}
\par}

\def\endproof{\null\hfill\endmark\endtrivlist}

\DeclareMathOperator{\Div }{div}

\def\pa{\partial}
\def\cal{\mathcal}
\let\mib=\boldsymbol

\def\R{\mathbb{R}}
\def\N{\mathbb{N}}

\def\eps{\varepsilon}

\def\mx{{\bf x}}
\def\mz{{\bf z}}
\def\mxi{{\mib \xi}}

\def\my{{\bf y}}
\def\mz{{\bf z}}

\def\ph{\varphi}
\def\Rd{{\mathbb{R}^{d}}}

\def\Sdmj{{\rm S}^{d-1}}
\def\Sd{{\rm S}^{d}}

\def\la{\lambda}

\def\mff{{\mathfrak f}}

\begin{document}
\title{Strong traces to degenerate parabolic equations}

\author{M.~Erceg}\address{Marko Erceg,
	Department of Mathematics, Faculty of Science, University of Zagreb, Bijeni\v{c}ka cesta 30,
	10000 Zagreb, Croatia}\email{maerceg@math.hr}

\author{D.~Mitrovi\'c}
\address{Darko Mitrovi\'c,
Faculty of Mathematics, University of
Vienna, Oscar Morgenstern platz 1,
1090 Vienna, Austria}\email{darko.mitrovic@univie.ac.at}




\begin{abstract}
We prove existence of strong traces at $t=0$ for quasi-solutions to (multidimensional) degenerate parabolic equations with no non-degeneracy conditions. In order to solve the problem, we combine the blow up method and a strong precompactness result for quasi-solutions to degenerate parabolic equations with the induction argument with respect to the space dimension.
\end{abstract}
\maketitle

\section{Introduction}

In the current contribution, we consider the advection diffusion
equation:

\begin{equation}
\label{d-p} \pa_t u+\Div_{\mx} \mff(u)=D^2_\mx \cdot
A(u) \,,
\end{equation} where $\mff\in C^1(\R; \R^d)$ and $A\in C^{1}(\R;\R^{d\times d})$  is a matrix symmetric at every point. Here we use the notation $D^2_\mx \cdot
A(u)=\sum\limits_{k,j}\pa^2_{x_k x_j}A_{kj}(u)$.
Usually, the equation above is written in 
the divergence form (non convenient for us at the moment)
\begin{equation*}
\pa_t u+\Div_{\mx} \mff(u)=\Div_{\mx}(a(u) \nabla u) \,,
\end{equation*} where $a(\lambda)=A'(\lambda)$ and $\Div_{\mx}(a(u) \nabla u)=\sum\limits_{k,j}^d \pa_{x_j}(a_{kj}(u)\pa_{x_k}u)$. 
The given equation is very common in applications as it describes phenomena
containing the combined effects of nonlinear convection, degenerate
diffusion, and nonlinear reaction. More precisely, the equation
describes a flow governed by

\begin{itemize}

\item the convection effects (bulk motion of particles) which are
represented by the first order terms;

\item diffusion effects which are represented by the second order
term and the matrix
$A'(\lambda)=[a_{ij}(\lambda)]_{i,j=1,\dots,d}$ describes direction and intensity of the diffusion.

\end{itemize} The equation is degenerate in the sense that the matrix $a(\lambda)=A'(\lambda)$ can be equal to zero in some directions, which are allowed to depend on $\lambda$. 
Roughly speaking, if this is the case, i.e.~if for some vector $\mxi
\in \R^{d}$ we have $\langle A'(\lambda)\mxi\,|\,\mxi\rangle =0$,
then diffusion effects do not exist for the state
$\lambda$ in the direction $\mxi$.

The equation appears in a broad spectrum of applications, such as e.g.~flow in porous media \cite{12}, sedimentation-consolidation processes \cite{6} and many others which we omit here (see the Introduction of \cite{karchen} for more details). 

Existence and uniqueness for the Cauchy problem corresponding to \eqref{d-p} is well established in quite general situations \cite{Car, CP, karchen}. 
The question of existence of strong traces for entropy solutions to \eqref{d-p} is however still open.
Before precisely formulating the problem and our results, let us first introduce the notion of 
quasi-solutions to \eqref{d-p} which is considered in this paper (see \cite{HKMP, NPS18, pan_jhde}).

By ${\rm sgn}$ we denote the signum function, while by $\mathcal{M}(\R^{d+1}_+)$ we denote the space of unbounded Radon measures $\nu$ on $\R^{d+1}_+$
	which are locally finite up to the boundary $t=0$, 
	i.e.~
$$	
\text{for any $T>0$ and compact set $K\subseteq \Rd$ it holds 
	$\operatorname{Var}\nu((0,T]\times K))<\infty $}\,,
$$ where $\operatorname{Var}\nu$ is the total
	variation of $\nu$.

\begin{definition}\label{def:quasisol}
	A measurable function $u$ defined on $\R^{d+1}_+:=(0,+\infty)\times\Rd$ is called a
	quasi-solution to \eqref{d-p} if 
	$u\in L^1_{loc}(\R^{d+1}_+)$,
	$\mff(u(\cdot))\in L^1_{loc}(\R^{d+1}_+;\Rd)$,
	$A(u(\cdot))\in L^1_{loc}(\R^{d+1}_+;\R^{d\times d})$, 
	and for any $\lambda\in \R$ the Kruzhkov type entropy equality holds
	\begin{align}
	\label{e-c} 
	\pa_t |u-\lambda| &+\Div_\mx\Bigl({\rm sgn}(u-\lambda)\bigl(\mff(u)-\mff(\lambda)\bigr)\Bigr)\\
	&-D^2_\mx\cdot\Bigl( {\rm sgn}(u-\lambda) \bigl(A(u)-A(\lambda)\bigr)\bigr) 
	= -\gamma(t,\mx,\lambda) \,, \nonumber
	\end{align} where $\gamma\in C(\R_\lambda;w\star-{\cal
		M}(\R_+^{d+1}))$ we call the quasi-entropy defect measure.
	 
\end{definition}

\begin{remark}
	We shall often consider $\gamma\in C(\R_\lambda;w\star-{\cal
		M}(\R_+^{d+1}))$ as a measure in $\lambda$ as well.
	More precisely, when we write $d\gamma(t,\mx,\lambda)$ in fact we 
	would think of $d\gamma(\lambda)(t,\mx)d\lambda$, where $d\lambda$ 
	is the Lebesgue measure. 
	Since $\gamma$ is continuous in $\lambda$, measure
	$\gamma\in {\cal M}(\R_+^{d+1}\times\R)$ is 
	locally finite up to the boundary $t=0$ as well. 
	Indeed, for any compact set $K_2\subseteq\R$, set $\{\gamma(\lambda) : \lambda\in K_2\}$
	is vaguely bounded in ${\cal M}(\R_+^{d+1})$, i.e.~with respect to the weak-$\star$
	topology.
	Thus, for any $T>0$ and any compact set $K_1\subseteq\Rd$ there exists 
	$c_{T,{K_1}}>0$ such that $\operatorname{Var}\bigl(\gamma(\lambda)\bigr)\leq c_{T,{K_1}}$
	(see e.g.~\cite[(13.4.2)]{JDieud}).
	Now, $\operatorname{Var}\gamma \leq c_{T,{K_1}} \operatorname{meas}(K_2)$,
	where $\operatorname{meas}(K_2)$ stands for the Lebesgue measure 
	of $K_2\subseteq\R$.
\end{remark}

The notion of quasi-solutions is introduced in \cite{pan_jhde} and it is a
generalization of the Kruzhkov-type admissibility concept (see e.g. \cite{Car, CP, Kru}). 
In a special
situation, i.e.~$\gamma\geq 0$, the quasi-solution is an entropy admissible solution that
singles out a physically relevant solution to the equation
\eqref{d-p} (see e.g. \cite{CP}).

Recently, several existence results of quasi-solutions to \eqref{d-p} in the case of irregular and heterogeneous fluxes were obtained (see e.g.~\cite{EMM20, HKMP, pan_arma}). All these results require 
suitable non-degeneracy conditions to be fulfilled (see \eqref{optimal}).
Such kind of assumptions are standard in the theory of velocity averaging lemmas \cite{EMM20,LM2, PS, TT} (see in particular \cite[(2.18)-(2.19)]{TT}), which is substantially used in the frame of the blow up method \cite{vass}. 
In the current contribution 
we shall rely on similar results from \cite{HKMP}.
However, let us emphasise that we are able to omit these non-degeneracy 
conditions in the final and main result of the paper (see Theorem 
\ref{main-traces}). 

Let us now recall that a function $u=u(t,\mx)$ has the (strong) trace $u_0=u_0(\mx)$ at $t=0$ if $L^1_{loc}-\lim\limits_{t\to 0} u(t,\cdot)=u_0$. 
 More precisely, we shall use the following definition.
\begin{definition}\label{traces}
Let $u\in L^1_{loc}(\R^{d+1}_+)$. A locally integrable function $u_0$ defined on $\R^{d}$ is called the strong trace of $u$ at $t=0$ if $\operatorname{ess\,lim}_{t\to 0^+}u(t,\cdot)=u_0$ 
in $L^1_{loc}(\Rd)$, i.e.~for some set $E\subseteq (0,\infty)$ of full Lebesgue 
measure and any relatively compact set $K\subset \subset \R^{d}$ it holds
\begin{equation}
\label{trf}
\lim\limits_{E\ni t\to 0} \|u(t,\cdot)-u_0\|_{L^1(K)}=0 \,.
\end{equation}

\end{definition} 

The strong traces appeared in the context of limit of hyperbolic relaxation toward scalar conservation laws \cite{natalini, tzavaras}. In particular, they appeared to be very useful related to the uniqueness of solutions to scalar conservation laws with discontinuous fluxes (see very restrictive list \cite{AKR, AM, Crasta} and references therein).

It is well known that the notion of weak solutions to \eqref{d-p} is 
not sufficient to ensure the existence of strong traces (see \cite{pan_jhdeB}).
Therefore, an
entropy requirement (at least in the weak form \eqref{e-c})
has to be used to establish such
a result.

One of the first results concerning the existence of strong traces was proved in \cite{chen} for the one dimensional situation. In the multidimensional case,  the result was obtained in \cite{vass} for entropy solutions to scalar conservation laws \cite{Kru} under non-degeneracy conditions. There, the basic technique for the proof -- the blow up method -- was introduced in this context. The results are further extended in \cite{pan_jhde, pan_jhdeB} for quasi-solutions to scalar conservation laws without the non-degeneracy conditions by combining the blow up techniques, H-measures, and induction with respect to the space dimension -- the technique that enabled avoiding the non-degeneracy conditions. In the frame of the technique, one introduces a change of variables which locally removes derivative with respect to one of the variables (i.e.~one of the flux components locally becomes zero and the corresponding variable becomes a parameter). Here, we are able to adjust and apply this methodology, therefore, we do not need non-degeneracy conditions eventually.

Let us remark in passing that existence of strong traces for entropy solutions for general multi-dimensional scalar conservation laws is still open in the sense that it is yet not known whether there exists $u_0\in L^\infty(\R^d)$ such that for an entropy solution to 
$$
\pa_t u+\Div_\mx \mff(t,\mx,u)=0 \;, \ \ \mff \in C^1(\R^+\times \R^d \times \R) \;,
$$ relation \eqref{trf} holds. For some more recent results we refer to \cite{NPS18}, where in particular the authors showed that some-kind 
of non-degeneracy conditions are necessary for the existence of 
strong traces in the case of rough fluxes.

As for the (entropy) solutions to degenerate parabolic equations 
(see e.g. \cite{Car, karchen, CP, VH}),
the first result for degenerate parabolic equations was obtained in \cite{Kwon09}, where the 
authors studied the case of scalar matrix $a\in C^1(\R)$. This essentially implies that the equation is locally either parabolic (if $a(\lambda)>0$) or hyperbolic (if $a(\lambda)=0$) which avoids the problem of dependence on $\lambda$ of directions in which the matrix $a$ degenerates. 
We avoided this problem in \cite{AM_jhde} by assuming that the matrix $a$ degenerates in fixed directions. Such matrices are called ultra-parabolic. We note that results from \cite{AM_jhde, Kwon09} do not cover e.g.~the matrix $a(\lambda)=\begin{bmatrix} \lambda^2 & \lambda\\
\lambda & 1 \end{bmatrix}$.
Recently, some progress is made in \cite{EMM20} where precisely the case 
of aforementioned matrix is resolved, while a general situation 
remained open.

An obvious problem with (general) degenerate parabolic equations
is inadequacy of the standard blow-up technique which involves scalings of the variables. Namely, if we are in the hyperbolic setting we use the scaling $(t,\mx) \mapsto (\eps t, \eps \mx)$ (the same with respect to both variables) \cite{vass}, while in the (ultra) parabolic setting, we need $(t,\mx) \mapsto (\eps t, \sqrt{\eps} \bar{\mx}+\eps \hat{\mx})$, $\mx=(\bar{\mx},\hat{\mx})$ \cite{AM_jhde}. This clearly causes problems if the equation changes type, since we cannot use the adapted scaling as in the ultra-parabolic case. In the current contribution, the method of the proof enables us to overcome these difficulties, thus obtaining the result for degenerate parabolic equations with diffusion matrix changing the directions of degeneracy. 

More precisely, we prove the following theorem. 

\begin{theorem}\label{main-traces} 
	Let $\mff\in C^1(\R;\R^{d})$ and  let $A\in C^{1}(\R;\R^{d\times d})$ be such that 
	 $a(\lambda):=A'(\lambda)$ is a positive semi-definite matrix.
	
	Then any bounded quasi-solution $u\in L^\infty(\R^{d+1}_+)$ to \eqref{d-p} admits the bounded strong trace at $t=0$, i.e.~there exists $u_0\in L^\infty(\R^{d})$ such that $$\operatorname{ess\,lim}_{t\to 0^+}u(t,\cdot)=u_0$$
	strongly in $L^1_{loc}(\Rd)$.
\end{theorem}

Using the truncation argument (see e.g.~\cite[Theorem 28]{HKMP}), we have the following extension of the previous theorem. 
\begin{corollary}\label{main-traces-cor} 
	Under assumptions of the previous theorem, let 
	$u\in L^\infty_{loc}(\R^+;L^p_{loc}(\R^{d}))$, for some $p>1$, 
	be a quasi-solution to \eqref{d-p}.
	Then $u$ admits the strong trace at $t=0$.
\end{corollary}

Let us now introduce the non-degeneracy conditions that we shall need in the sequel.

\begin{definition}\label{def:non-deg-ab}
We say that the flux $\mff$ and the diffusion matrix $A$ satisfy \emph{non-degeneracy conditions} on $(a,b)\subseteq\R$ if there exists no interval $(\alpha, \beta)\subseteq (a,b)$, $\alpha<\beta$, such that for some $(\xi_0,\mxi) \in \Sd$ it holds
\begin{equation}
\label{optimal}
 \xi_0+\langle \mff'(\lambda)\,|\,\mxi \rangle \,=\, \langle A'(\lambda) \mxi\,|\,\mxi \rangle \,=\,0 \;, \ \ \lambda \in (\alpha, \beta) \;,
\end{equation}
i.e.~for any $(\xi_0,\mxi) \in \Sd$ functions
$\lambda\mapsto \xi_0+\langle \mff'(\lambda)\,|\,\mxi \rangle$,
$\lambda\mapsto \langle A'(\lambda) \mxi\,|\,\mxi \rangle$ are not 
simultaneously equal to zero on non-degenerate intervals contained 
in $(a,b)$. 
\end{definition}

In \eqref{optimal}, as well as in the rest of the paper, we use $\langle\cdot\,|\,\cdot\rangle$
to denote scalar product on $\Rd$ and $\mathrm{S}^d$ for the unit sphere in $\R^{d+1}$.

\smallskip

Let us now briefly describe the solution procedure. 
\begin{enumerate}[label=(\alph*) \,]
\item We rewrite the equation in the kinetic formulation and prove existence of weak traces for the solution of the kinetic equation 
(see Theorem \ref{w-t}).

\item We reduce the problem to the problem of existence of strong traces 
for $v:=s_{\alpha,\beta}(u)$, $\alpha<\beta$, where 
$s_{\alpha,\beta}(u)=\max\{ \alpha, \min\{u,\beta\}\}$
(see steps II--IV of the proof of Theorem \ref{main-traces}). 

\item We prove existence of strong traces for $v$ under 
\eqref{optimal} and for matrices $a:=A'$ of the form 
\begin{equation}
\label{wantedform}
a(\lambda)=\begin{bmatrix} \tilde{a}(\lambda) & 0\\
0 & 0
\end{bmatrix} \,,
\end{equation}
where $\tilde a(\lambda)\in\R^k\times\R^k$, for some $k\in\{0,1,\dots,d\}$, 
is such that for any $\mxi\in\R^k\setminus\{0\}$ the mapping
$\lambda \mapsto \langle \tilde a(\lambda)\mxi\,|\,\mxi \rangle$
is not identically equal to zero on non-degenerate intervals
(see Theorem \ref{sc-ndeg}).

\item We remove the non-degeneracy assumption on the flux in the previous 
step by reducing the dimension of the equation (see Theorem \ref{T2}). 

\item If the matrix $a$ does not satisfy assumptions given in (c) for 
$k=d$, then there exists $\mxi^*\in \Sdmj$ and an interval $(\alpha,\beta)$ such that  
\begin{equation*}
\langle a(\lambda)\mxi^*\,|\,\mxi^* \rangle=0 \;, \quad \lambda\in (\alpha,\beta) \,.
\end{equation*} 
This means that the matrix $a$ degenerates in the direction 
$\mxi^*$ uniformly for $\lambda\in (\alpha,\beta)$. 
By introducing appropriate change of variables, 
we can assume that $\mxi^*=e_d=(0,\dots,0,1)$. 
From here and since the matrix $a$ is non-negative definite, 
we have $a_{jd}(\lambda)=a_{dj}(\lambda)=0$, $j=1,2,\dots,d$, 
$\lambda \in (\alpha,\beta)$. Thus, the matrix $a$ has the form 
\eqref{wantedform} for $k=d-1$. If $\tilde a$ does not satisfy 
the non-degeneracy assumption from (c), we continue with this reduction
(which will finish in finite number of steps).
Finally, we apply (d) (see step I of the proof of Theorem \ref{main-traces}).  
\end{enumerate}

In Section 2 we prove existence of weak traces and develop a sufficient 
condition for the existence of strong traces in terms of rescaled sequences.
In Section 3 existence of strong traces is obtained in the case when 
the diffusion matrix is of the block structure \eqref{wantedform} and 
satisfies non-degeneracy conditions (point (d) in the above plan). 
The main result, existence of traces in the general situation without non-degeneracy conditions (Theorem \ref{main-traces}), 
is proved in Section 4. 
The paper is closed by few concluding remarks in the final section.

\section{Weak traces}

In this section we prove that any quasi-solution to 
\eqref{d-p} admits a weak trace. 
More precisely, we shall first develop the result for bounded 
quasi-solutions, while the general case is treated 
by truncating solutions using
\begin{equation}
\label{ab}
s_{\alpha,\beta}(u)=\max\{ \alpha, \min\{u,\beta\}\} \;.
\end{equation}
Thus, it is sufficient to prove that the weak trace is in fact the 
strong trace. We close this section with one such condition 
in terms of the strong convergence of a rescaled (blow-up) sequence.

In this paper we chose to use the kinetic reformulation (cf.~\cite{CP, LPT})
of the quasi-solution concept, which is not essential 
(see Concluding remarks).

\begin{theorem}
\label{T7}
Denote $f:=\mff'$ and $a:=A'$.
If function $u$ is a bounded quasi-solution to \eqref{d-p}, then
the function
\begin{equation}
\label{equil} 
h(t,\mx,\lambda)={\rm sgn}\bigl(u(t,\mx)-\lambda\bigr)=
 	-\pa_\lambda |u(t,\mx)-\lambda|
\end{equation} is a weak solution to the following linear equation:
\begin{align}
\label{k-1} \pa_t h + \Div_\mx \bigl( f(\lambda) \,h\bigr)
	-D^2_\mx\cdot\bigl( a(\lambda) \,h\bigr) = \pa_\lambda 
	\gamma(t,\mx,\lambda) \,.
\end{align} Conversely, if \eqref{equil} (with bounded $u$) satisfies \eqref{k-1}, then $u$ is the quasi-solution to \eqref{d-p}.
\end{theorem}

\begin{proof}
It is enough to find derivative of \eqref{e-c} with respect to $\lambda\in
\R$ to obtain \eqref{k-1}. Conversely, it is sufficient to integrate \eqref{k-1} on the interval $(-\|u\|_{\infty},\lambda)$ and to take into account that $\gamma\in C(\R;{\cal M}(\R^{d+1}_+))$ to reach to \eqref{e-c}. Details of the procedure can be found in e.g.~\cite{CP} in a slightly less general setting. 
\end{proof}

We have the following theorem.

\begin{theorem}
\label{w-t} 
Let $h\in L^\infty(\R^{d+1}_+\times \R)$ be a
distributional solution to \eqref{k-1} and let us define
\begin{equation}\label{eq:set_E}
\begin{aligned}
E=\Bigl\{t\in \R^+: \; (t,\mx,\lambda)& \text{  is a
Lebesgue point of}\\ & h(t,\mx,\lambda) \text{  for a.e.} \, (\mx,\lambda) \in \R^{d}\times \R \Bigr\} \,.
\end{aligned}
\end{equation}

Then there exists $h_0\in L^\infty(\R^{d+1})$, such that
$$
h(t,\cdot, \cdot) \xrightharpoonup{\;*\;} h_0\,,  \mbox{ \rm weakly-$\star$ in } L^\infty(\R^{d+1}) \,,
	\mbox{ \rm as } t\to 0\,,\; t\in E\,.
	$$
\end{theorem}

\begin{proof}Note first that $E$ is of full measure. Since $h\in L^\infty(\R^{d+1}_+\times \R)$, the family $\{
h(t,\cdot,\cdot)\}_{t\in E}$ is bounded in $L^\infty(\R^{d+1})$.
             Due to the weak-$\star$ precompactness of $L^\infty(\R^{d+1})$, there exists a sequence $\{t_m\}_{m\in \N}$ in $E$ such that $t_m\to 0$ as $m\to \infty$,
             and $h_0\in L^\infty(\R^{d+1})$, such that
                         \begin{equation}\label{wc}
                                  h(t_m, \cdot,\cdot) \xrightharpoonup{\;*\;} h_0 \,, \mbox{ weakly-$\star$ in $L^\infty(\R^{d+1})$, as } m\to \infty \,.
                         \end{equation}
             For $\phi\in C_c^\infty(\R^{d})$, $\rho\in C^1_c(\R)$, denote
             $$
               I(t):=\int_{\R^{d+1}}h(t,\mx,\lambda)\rho(\lambda)\phi(\mx)\,d\mx d\lambda \,, \quad t\in E\,.
             $$
             With this notation, \eqref{wc} means that
                         \begin{equation}\label{wc1}
                                \lim_{m\to \infty} I(t_m)= \int_{\R^{d+1}} h_0(\mx,\lambda)\rho(\lambda)\phi(\mx)\, d\mx d\lambda = :I(0) \,.
                         \end{equation}
             Now, fix $\tau\in E$ and notice that for the regularization $I_\eps=I\star \omega_\eps$, where $\omega_\eps$ is the standard convolution kernel, it holds
             $$
             \lim\limits_{\eps \to 0} I_\eps(\tau)=I(\tau).
             $$
             
Then, fix $m_0\in \N$, such that $E\ni t_m\leq \tau$, for $m\geq m_0$, and remark that
 \begin{equation}
 \begin{split}\nonumber
 I(\tau)- I(t_m) &= \lim\limits_{\eps\to 0}\int_{t_m}^\tau I_\eps'(t)\,dt 
 \\
 &= \sum_{j=1}^d \int_{(t_m,\tau] \times \R^{d+1}} h(t,\mx,\lambda)f_j(t,\mx,\lambda) \rho(\lambda)\pa_{x_j}\phi(\mx)\, dt\, d\mx \,d\lambda \\
 &\qquad + \sum_{j,k=1}^d \int_{(t_m,\tau] \times \R^{d+1}} h(t,\mx,\lambda) a_{jk}(\lambda) \rho(\lambda) \pa_{x_j x_k}\phi(\mx)\, dt\, d\mx \,d\lambda \\
 &\qquad - \int_{(t_m,\tau] \times \R^{d+1}} \phi(\mx) \rho'(\lambda)\, d\gamma(t,\mx,\la) \,,
 \end{split}
 \end{equation}
 where we have used that $h$ is a distributional solution to \eqref{k-1}.
Hence, passing to the limit as $m\to\infty$, and having in mind \eqref{wc1} and the fact that
$\gamma$ is locally finite up to the boundary $t=0$, we obtain
  \begin{equation}\begin{split}\nonumber
  I(\tau)- I(0) &= \sum_{j=1}^d \int_{(0,\tau] \times \R^{d+1}} h(t,\mx,\lambda) f_j(t,\mx,\lambda) \rho(\lambda)\pa_{x_j}\phi(\mx) \, dt\, d\mx \,d\lambda \\
  &\qquad + \sum_{j,k=1}^d \int_{(0,\tau] \times \R^{d+1}} h(t,\mx,\lambda) a_{jk}(\lambda) \rho(\lambda)  \pa_{x_j x_k}\phi(\mx)\, dt\, d\mx \,d\lambda \\
  &\qquad - \int_{(0,\tau] \times \R^{d+1}} \phi(\mx) \rho'(\lambda) \, d\gamma (t,\mx,\la) \,.
  \end{split}\end{equation}
  The right hand side clearly tends to zero as $\tau\to 0$.
             Thus, for all $\phi\in C_c^\infty(\R^{d+1})$ and $\rho\in C^1_c(\R)$
             we have $\lim_{E\ni\tau\to 0}I(\tau)=I(0)$, i.e.
             $$
             \lim_{E\ni\tau\to 0 } \int_{\R^{d+1}} h(\tau,\mx,\lambda)\rho(\lambda) \phi(\mx)\, d\mx d\lambda = \int_{\R^{d+1}} h_0(\mx,\lambda)\rho(\lambda) \phi(\mx) \, d\mx d\la \,.
             $$
Having in mind that $h(\tau, \cdot,\cdot)$, $\tau\in E$, 
is bounded, and that $C_c^\infty(\R^{d+1})$ is dense in
$L^1(\R^{d+1})$, we complete the proof.
\end{proof}

\begin{remark}\label{rem:w-t_u}
If $u$ is a bounded quasi-solution to \eqref{d-p}, 
by the previous theorem it admits the weak trace
(cf.~\cite[Corollary 2.2]{pan_jhde} for the case $A=0$).

Indeed, let $M>0$ be such that 
$$
|u(t,\mx)| \leq M \;, \quad \hbox{a.e} \ (t,\mx)\in\R^{d+1}_+ \;.
$$
Then, by the definition of $h$ (it is the sign function; see \eqref{equil})
$$
\int_{-M}^{M} h(t,\mx,\lambda)\,d\lambda 
	= \int_{-M}^{M} {\rm sgn}\bigl(u(t,\mx)-\lambda\bigr)\,d\lambda
	= \int_{-M}^{u(t,\mx)} d\lambda - \int_{u(t,\mx)}^{M} d\lambda
	= 2u(t,\mx) \,.
$$
Thus, the claim follows by Theorem \ref{w-t} by noting that the characteristic function 
$\chi_{[-M,M]}$ of the interval $[-M,M]$ is in $L^1(\R)$.
More precisely, we have
\begin{equation}\label{eq:u_w-t}
u(t,\cdot) \xrightharpoonup{\;*\;} u_0:=\frac{1}{2}\int_{-M}^{M} h_0(\cdot,\lambda)\,d\lambda
\end{equation}
weakly-$\star$ in $L^\infty(\R^{d})$ as $t\to 0$, $t\in E$.
\end{remark}

\begin{remark}\label{rem:w-t_u_B}
In a similar fashion as in the previous remark, we have for any $f\in L^1(\R)$
$$
\int_{-M}^{M} f(\lambda) h(t,\mx,\lambda)\,d\lambda \xrightharpoonup{\;*\;} \int_{-M}^{M} f(\lambda) h_0(\mx,\lambda)\,d\lambda
$$ 
weakly-$\star$ in $L^\infty(\R^{d})$ as $t\to 0$, $t\in E$.

In particular, for $\lambda\mapsto\lambda\chi_{[-M,M]}(\lambda)\in L^1(\R)$ we have
$$
\int_{-M}^{M} \lambda h(t,\mx,\lambda)\,d\lambda 
	= \int_{-M}^{u(t,\mx)} \lambda \,d\lambda - \int_{u(t,\mx)}^{M} \lambda\,d\lambda
	= u(t,\mx)^2 -M^2 \,
$$ and thus
\begin{equation}\label{eq:uu_w-t}
u^2(t,\cdot) \xrightharpoonup{\;*\;} u_1:=\int_{-M}^{M} \lambda h_0(\cdot,\lambda)
	\,d\lambda +M^2
\end{equation} weakly-$\star$ in $L^\infty(\R^{d})$ as $t\to 0$, $t\in E$.

If one can get that $u_1=u_0^2$, the strong convergence 
(i.e.~the strong trace) can be obtained from the above weak convergences. 
Namely, for an arbitrary $\varphi\in C_c(\Rd)$ by \eqref{eq:u_w-t}--\eqref{eq:uu_w-t}
we have
\begin{equation}\label{eq:u_strong_conv}
\begin{aligned}
\lim_{E\ni t\to 0}\int_\Rd \Bigl(u(t,\mx)- & u_0(\mx)\Bigr)^2 \varphi(\mx) \,d\mx \\
&= \lim_{E\ni t\to 0}\int_\Rd \Bigl(u(t,\mx)^2-2u(t,\mx)u_0(\mx)+u_0(\mx)^2\Bigr)\varphi(\mx) 
\,d\mx \\
&= \int_\Rd \Bigl(u_1(\mx)-u_0(\mx)^2\Bigr)\varphi(\mx) \,d\mx \,.
\end{aligned}
\end{equation}

An obvious way to show identity $u_1=u_0^2$ is to prove that 
for a sequence $(t_m)$ in $E$ converging to zero the sequence 
$u_m:=u(t_m,\cdot\,)$ converges strongly in $L^1_{loc}(\Rd)$. 
However, it will be useful to develop another sufficient condition
in terms of the strong convergence of certain rescaled (sub)sequences,
which is more adequate to the equation \eqref{d-p}
(see Lemma \ref{equivalence} below).
\end{remark}

Now we use the rescaling procedure (or the so-called blow-up method) 
in order to 
obtain a sufficient condition for the existence of the strong trace. 
More precisely, for an arbitrary (but fixed) $k\in\{0,1,\dots,d\}$ 
let us consider the sequence
\begin{equation} 
\label{scal}      
u\Bigl(\frac{\hat  t}{m}, \frac{\tilde \mx}{\sqrt{m}}+\tilde{\my},\frac{\bar{\mx}}{{m}}+\bar{\my} \Bigr)	\,, \ \ m\in \N \,,
\end{equation}
where $\tilde\mx, \tilde \my\in \R^k$ and 
$\bar{\mx},\bar{\my}\in\R^{d-k}$.
We shall use a shorthand $\hat{\mx}=(\tilde\mx,\bar{\mx})$ and 
$\my=(\tilde{\my},\bar{\my})$. 
The following lemma holds.

\begin{lemma}
\label{equivalence}
Let $u$ be a bounded quasi-solution to \eqref{d-p}.
Assume that the sequence given by \eqref{scal} converges toward 
$(\hat{t},\hat{\mx},\my)\mapsto u_0(\my)$ 
in $L^1_{loc}(\R^{d+1}_+\times \R^d)$ along a subsequence, 
where $u_0$ is the weak trace of $u$ (see Theorem \ref{w-t}).
Then function $u$ admits the strong trace at $t=0$. 

Conversely, if a quasi-solution $u$ admits the strong trace 
at $t=0$, then the sequence given by \eqref{scal} converges 
toward the strong trace.
\end{lemma}

Before proving the lemma, let us explain how the strong 
convergence of \eqref{scal} can be obtained. 
For a fixed $\my=(\tilde{\my},\bar{\my})\in \Rd$, let us 
change the variables in \eqref{e-c} in the following way:
\begin{equation}
\label{blow-up}
(t,\mx)= \Bigl(\frac{\hat t}{m}, \frac{\tilde \mx}{\sqrt{m}}+\tilde{\my},\frac{\bar{\mx}}{{m}}+\bar{\my}\Bigr) \,.
\end{equation} 
If $u$ is a quasi-solution to \eqref{d-p}, then the rescaled function
\eqref{scal}, denoted here by $u_m$, satisfies (we purposely label the two rows below)
\begin{align}
\pa_{\hat{t}} |u_m -\lambda|
	&+\Div_{\bar\mx}\Bigl({\rm sgn}(u_m-\lambda)
	\bigl(\bar\mff(u_m)-\bar\mff(\lambda)\bigr)\Bigr)
	-\sum\limits_{s,j=1}^k \pa_{\tilde{x}_j \tilde{x}_s}
	\Bigl({\rm sgn}(u_m-\lambda) 
	\bigl({A}_{sj}(u_m)-{A}_{sj}(\lambda)\bigr)\Bigr)
	\nonumber \\
&-2\sqrt{m}\sum\limits_{j=1}^{k}\sum\limits_{s=1}^{d-k} 
	\pa_{\tilde{x}_j \bar{x}_s}\Bigl(
	{\rm sgn}(u_m-\lambda) 
	\bigl({A}_{(k+s)j}(u_m)-{A}_{(k+s)j}
	(\lambda)\bigr)\Bigr) \label{eq:unbdd_terms_non-diag} \\
&-m\sum\limits_{s,j=1}^{d-k} \pa_{\bar{x}_j \bar{x}_s}\Bigl(
	{\rm sgn}(u_m-\lambda) 
	\bigl({A}_{(k+s)(k+j)}(u_m)-{A}_{(k+s)(k+j)}
	(\lambda)\bigr)\Bigr) \label{eq:unbdd_terms_diag} \\
&= -\frac{1}{m}\gamma_m^\my(\hat{t},\hat{\mx},\lambda) 
	-\frac{1}{\sqrt{m}}\Div_{\tilde\mx}\Bigl({\rm sgn}(u_m-\lambda)
	\bigl(\tilde\mff(u_m)-\tilde\mff(\lambda)\bigr)\Bigr)\,,
	\nonumber
\end{align}  
where $\mff=(\,\tilde{\mff},\bar{\mff}\,)$ and the equality between 
$\gamma$ and $\gamma_m^\my$ is understood in the sense of distributions:
\begin{equation}
\label{**}
\langle \gamma_m^\my,\varphi \rangle=m^{k/2+d-k+1} 
\int_{\R^{d+1}_+\times\R}\varphi\bigl(m\, t,\sqrt{m} (\tilde{\mx}-\tilde{\my}),{m} (\bar{\mx}-\bar{\my}),\lambda\bigr)d\gamma({t},{\mx},\lambda) \,.
\end{equation}
Now the strong convergence of $(u_m)$ is obtained from the equation
above by applying a compactness result from \cite{HKMP},
which is done in the following section. 
Let us just remark that in order to do so, 
we need to choose $k$ such that unbounded terms 
\eqref{eq:unbdd_terms_non-diag}--\eqref{eq:unbdd_terms_diag} do not appear
and that the remaining left hand side satisfies non-degeneracy conditions.
\smallskip

\noindent \textbf{Proof of Lemma \ref{equivalence}:}
By the assumptions we have for any non-negative $\varphi\in C_c(\R^{d+1}_+\times \R^{d})$ along the subsequence from the formulation of the lemma
\begin{align*}
\lim\limits_{m\to \infty}\int_{\Rd\times\R^{d+1}_+}
	\varphi(\hat{t},\hat{\mx},\my)
	\biggl| u\Bigl(\frac{\hat  t}{m}, \frac{\tilde \mx}{\sqrt{m}}+\tilde{\my},\frac{\bar{\mx}}{{m}}+\bar{\my} \Bigr)-u_0(\my) \biggr| \,d\my d\hat{\mx} d \hat{t}=0 \,,
\end{align*}
where $\hat{\mx}=(\tilde{\mx},\bar{\mx})$ and $\my=(\tilde{\my},\bar{\my})$.

Taking into account the change of variables $\tilde{\mz}=\frac{\tilde \mx}{\sqrt{m}}+\tilde{\my}$ and $\bar{\mz}=\frac{\bar \mx}{{m}}+\bar{\my}$
($\mz=(\tilde{\mz},\bar{\mz})$) with respect to $\my$, the previous limit reads 
\begin{equation}
\label{****}
\begin{split}
0 &=\lim\limits_{m\to\infty}\int_{\R^{d}\times\R^{d+1}_+}
	\varphi\Bigl(\hat{t},\hat{\mx},\tilde{\mz}-\frac{\tilde \mx}
	{\sqrt{m}},\bar{\mz}-\frac{\bar \mx}{{m}}\Bigr)
	\biggl|u\Bigl(\frac{\hat t}{m},\mz\Bigr)
	-u_0\Bigl(\tilde{\mz}-\frac{\tilde \mx}{\sqrt{m}},
	\bar{\mz}-\frac{\bar \mx}{{m}}\Bigr) 
	\biggr| \,d\mz d\hat{\mx} d\hat{t} \\
&= \lim\limits_{m\to\infty}\int_{\R^{d}\times\R^{d+1}_+}
	\varphi(\hat{t},\hat{\mx},\mz)
	\biggl|u\Bigl(\frac{\hat t}{m},\mz\Bigr)
	-u_0(\mz) 
	\biggr| \,d\mz d\hat{\mx} d\hat{t}\;.
\end{split}
\end{equation}

Therefore, due to arbitrariness of $\varphi$ in \eqref{****}, we conclude 
$$
u\Bigl(\frac{\hat{t}}{m},\mz\Bigr)\to u_0(\mz)\,, \ \ m\to \infty
$$ 
in $L^1_{loc}(\R^{d+1}_+)$ along the subsequence from the formulation of the lemma. This means that (for another subsequence not relabelled) 
there exists $\hat{E}\subseteq\R^+$ of full measure such that for any $\hat{t}\in\hat{E}$
we have
\begin{equation}
\label{sc1}
u\Bigl(\frac{\hat{t}}{m},\mz\Bigr)\to u_0(\mz), \ \ m\to \infty
\end{equation}  in $L^1_{loc}(\R^{d})$. 
It is easy to see that for $E$ given by \eqref{eq:set_E} set 
$E_\infty:=\bigcap\limits_{m\in\N} {m} E$ is of full measure
(since $E$ is of full measure). 
Thus, the intersection $\hat{E}\cap E_\infty$ is non-empty (in fact it is a set of full 
measure as well), so we can choose $\hat{t}\in\R^+$ such that \eqref{sc1} holds and 
$\frac{\hat{t}}{m}\in E$, $m\in\N$.

Now, choose $\rho(\lambda)=\lambda \,\chi_{[-M,M]}(\lambda)$ where $\chi_{[-M,M]}(\lambda)$ is the characteristic function of the interval $[-M,M]$. It holds according to Theorem \ref{w-t} (see also Remark \ref{rem:w-t_u_B})
$$
u^2(t,\mx) = \int_{-M}^M \lambda h(t,\mx,\lambda) \,d\lambda + M^2
	\xrightharpoonup{\;*\;} \int_{-M}^M \lambda h_0(\mx,\lambda) \,d\lambda +M^2 
	=: u_1(\mx) 
$$ 
in $L^\infty(\R^{d})$ as  $E\ni t\to 0$. 
Since the weak-$\star$ convergence in $L^\infty(\Rd)$ implies the weak convergence in $L^1_{loc}(\Rd)$, and since weak and strong limits coincide, 
from here and \eqref{sc1} we see that it must be $u_1=u_0^2$. 
Finally, by \eqref{eq:u_strong_conv} (see Remark \ref{rem:w-t_u_B})
we have
$$
u(t,\cdot) \rightarrow u_0 
$$
in $L^2_{loc}(\Rd)$ as $t\to 0$, $t\in E$, which implies the convergence in $L^1_{loc}$.
Hence, $u_0$ is the strong trace. 

Conversely, assume that $u_0$ is the strong trace to a quasi-solution $u$ of \eqref{d-p}. By Definition \ref{traces}, it holds
\begin{equation*}
\operatorname{ess\,lim}_{t\to 0^+} \|u(t,\cdot)-u_0\|_{L^1(K)} 
	\;, \ \ K\subset\subset \R^d \,.
\end{equation*} 
In particular, for almost every $t\in \R^+$ and every $\varphi\in C_c(\R^{d+1}_+\times \R^{d})$, it holds according to the Lebesgue dominated convergence theorem:
$$
\lim\limits_{m\to\infty}
	\int_{\R^{d}\times\R^{d}}\varphi\Bigl({t},\hat{\mx},
	\tilde{\mz}-\frac{\tilde \mx}{\sqrt{m}},
	\bar{\mz}-\frac{\bar \mx}{{m}}\Bigr) 
	\biggl|u\Bigl(\frac{t}{m},\mz\Bigr)-u_0(\mz)\biggr| 
	\,d\mz d\mx =0
$$  
(here $\mx=(\tilde{\mx},\bar{\mx})$, $\mz=(\tilde{\mz},\bar{\mz})$ as usual). 
Applying the change of variables here, we have
$$
\lim\limits_{m\to \infty}\int_{\R^{d}\times\R^{d}}
	\varphi({t},\hat{\mx},\my) 
	\biggl|u\Bigl(\frac{t}{m},\frac{\tilde{\mx}}{\sqrt{m}}+\tilde{\my},
	\frac{\bar{\mx}}{m}+\bar{\my}\Bigr)
	-u_0\Bigl(\frac{\tilde{\mx}}{\sqrt{m}}+\tilde{\my},
		\frac{\bar{\mx}}{m}+\bar{\my}\Bigr)\biggr| \,d\my d{\mx} =0 \,, 
$$
where $\my=(\hat{\my},\tilde{\my})$. Using the fact that almost every point of an $L^1$-function is the Lebesgue point, we have from the Lebesgue dominated convergence theorem
$$
\lim\limits_{m\to \infty}\int_{\R^{d}\times\R^{d}}\varphi({t},\hat{\mx},\my)
	\biggl|u\Bigl(\frac{t}{m},\frac{\tilde{\mx}}{\sqrt{m}}+\tilde{\my},
	\frac{\bar{\mx}}{m}+\bar{\my}\Bigr)	
	-u_0(\my)\biggr| \,d\my d{\mx} = 0 \,. 
$$
Integrating the latter over $[0,\infty)$ and using the Lebesgue dominated convergence theorem to exchange the limit and integral, we get
$$
\lim\limits_{m\to \infty}\int_{\R^{d}\times\R_+^{d+1}}
	\varphi({t},\hat{\mx},\my) 	\biggl|u\Bigl(\frac{t}{m},\frac{\tilde{\mx}}{\sqrt{m}}+\tilde{\my},
	\frac{\bar{\mx}}{m}+\bar{\my}\Bigr)
	-u_0(\my)\biggr| \, d\my d{\mx} dt = 0 \;. 
$$
Thus, from the arbitrariness of $\varphi$, we have the conclusion of the second part of the lemma.
\endproof

\begin{remark}\label{rem:sab-sol}
For a quasi-solution $u$ to \eqref{d-p} and any $\alpha,\beta\in\R$, 
$\alpha<\beta$, 
function $v:=s_{\alpha,\beta}(u)$, where $s_{\alpha,\beta}$ is given by \eqref{ab}, 
is a bounded quasi-solution to \eqref{d-p}.
Indeed, if $u$ satisfies \eqref{e-c} with measure $\gamma$, 
then $v$ satisfies the same identity with 
$$
\tilde\gamma (t,\mx,\lambda) := 
	\gamma(t,\mx,s_{\alpha,\beta}(\lambda))-\frac{1}{2}\bigl(\gamma(t,\mx,\alpha)
	+\gamma(t,\mx,\beta)\bigr) \;, 
$$
i.e.
\begin{equation}\label{eq:sab_e-c}
\begin{aligned}
	\pa_t |v-\lambda| &+\Div_\mx\Bigl({\rm sgn}(v-\lambda)\bigl(\mff(v)-\mff(\lambda)\bigr)\Bigr)\\
	&-D^2_\mx\cdot\Bigl( {\rm sgn}(v-\lambda) \bigl(A(v)-A(\lambda)\bigr)\bigr) 
	= -\tilde\gamma(t,\mx,\lambda)
\end{aligned}
\end{equation}
(cf.~\cite[Section 7]{pan_jhde}).
Furthermore, taking derivative with respect to $\lambda$ we 
obtain that 
$$
h^{\alpha,\beta} (t,\mx,\lambda)=\operatorname{sgn}\bigl(v(t,\mx)-\lambda\bigr)
$$
satisfies 
\begin{equation}\label{eq:sab_k-f}
\pa_t h^{\alpha,\beta} + \Div_\mx \bigl( f(\lambda) \,h^{\alpha,\beta}\bigr)
	-D^2_\mx\cdot\bigl( a(\lambda) \,h^{\alpha,\beta}\bigr) = \pa_\lambda
	\bigl(\gamma(t,\mx,s_{\alpha,\beta}(\lambda))\bigr) \;.
\end{equation}

Therefore, if a quasi-solution $u$ to \eqref{d-p} is not 
bounded, we can replace it with $s_{\alpha,\beta}(u)$ and apply 
all the results valid for bounded quasi-solutions. 
Since functions $w\mapsto s_{\alpha,\beta}(w)$ are continuous,
it is important to notice that the set $E$ given by \eqref{eq:set_E}
is the same if we replace $u$ by $s_{\alpha,\beta}(u)$ in its definition.

Furthermore, if $s_{\alpha,\beta}(u)$ admits the strong trace
denoted by $u_0^{\alpha,\beta}$, then for any 
$\alpha',\beta'\in\R$, $\alpha\leq \alpha'<\beta'\leq\beta$,
the function $s_{\alpha',\beta'}(u)$ admits the strong trace 
as well, and its strong trace is equal to 
$s_{\alpha',\beta'}(u_{\alpha,\beta}^0)$
(this easily follows by noting that the strong convergence implies 
the convergence a.e.~along a subsequence).
In particular, if $u$ admits the strong trace on intervals 
$(\alpha_1,\beta_1)$, $(\alpha_2,\beta_2)$, i.e.~functions
$s_{\alpha_1,\beta_1}(u)$ and $s_{\alpha_2,\beta_2}(u)$
admit the strong traces $u_0^{\alpha_1,\beta_2}$ and $u_0^{\alpha_2,\beta_2}$
respectively, then if 
$\alpha:=\max\{\alpha_1,\alpha_2\}<\min\{\beta_1,\beta_2\}:=\beta$,
we have 
$s_{\alpha,\beta}(u_0^{\alpha_1,\beta_2})=
s_{\alpha,\beta}(u_0^{\alpha_2,\beta_2})$.

Moreover, if we have existence of the strong trace only for 
$s_{\alpha',\beta'}(u)$, denoted by $u_0^{\alpha',\beta'}$, 
and if we denote by $u_0^{\alpha,\beta}$ the weak trace
of $s_{\alpha,\beta}(u)$, $(\alpha',\beta')\subseteq (\alpha,\beta)$,
then we can still conclude that $u_0^{\alpha',\beta'}=
s_{\alpha',\beta'}(u_0^{\alpha,\beta})$.
The proof of this statement is elaborated in Remark 
\ref{rem:weak-strong-trace}.
\end{remark}

\section{Existence of strong traces for quasi-solutions to \eqref{d-p} under conditions \eqref{wantedform}}

In this section we study \eqref{d-p} under additional 
assumptions on the structure of matrix $A$. 
Namely, throughout the section we assume (see also \eqref{wantedform}):
\begin{itemize}
\item[(A1)] There exists $k\in\{0,1,\dots,d\}$, 
$\alpha,\beta\in\R$, $\alpha<\beta$, and
$\tilde a\in C([\alpha,\beta];\mathbb{R}^{k\times k})$ such that 
\begin{equation}\label{eq:atilde-block}
a(\lambda)=A'(\lambda) 
	= \begin{bmatrix}\tilde{a}(\lambda)&0\\ 0&0\end{bmatrix}
	\;, \ \lambda\in [\alpha,\beta] \;,
\end{equation}
and $\tilde a$ satisfies the non-degeneracy condition on $(\alpha,\beta)$,
i.e.~there exists no interval $(\alpha_1,\beta_1)\subseteq (\alpha,\beta)$, 
$\alpha_1<\beta_1$, 
such that for some $\tilde\mxi\in\R^k\setminus\{0\}$ the mapping 
$\lambda\mapsto \langle \tilde a(\lambda)\tilde\mxi\,|\,\tilde\mxi\rangle$
is equal to zero on the whole interval $(\alpha_1,\beta_1)$. 
\end{itemize}

It is obvious that in the case $k=0$ (i.e.~$a=0$ on $[\alpha,\beta]$)
the assumption above is trivially satisfied. 

If $u$ is a quasi-solution to \eqref{d-p}, then for $\lambda\in \R$
by Remark \ref{rem:sab-sol} (see \eqref{eq:sab_k-f}) holds
\begin{equation}\label{eq:sab-ultrapar-k-f}
\pa_t h^{\alpha,\beta} + \Div_\mx \bigl( f(\lambda) \, h^{\alpha,\beta}\bigr)
	-D^2_{\tilde \mx}\cdot\bigl( \tilde a(\lambda) \, h^{\alpha,\beta}\bigr) 
	= \pa_\lambda
	\gamma(t,\mx,s_{\alpha,\beta}(\lambda)) \;,
\end{equation}
where $\mx=(\tilde{\mx},\bar{\mx})\in\R^k\times\R^{d-k}$ and 
$h^{\alpha,\beta}$ is associated to $v:=s_{\alpha,\beta}(u)$,
i.e.~$h^{\alpha,\beta}(t,\mx,\lambda)=\operatorname{sgn}\bigl(v(t,\mx)-\lambda\bigr)$.
Indeed, for $\lambda\in[\alpha,\beta]$ the claim follows by \eqref{eq:atilde-block}.
On the complement, $\lambda\not\in[\alpha,\beta]$, $h^{\alpha,\beta}$ is a 
constant function, hence \eqref{eq:sab-ultrapar-k-f} trivially holds since it
reads $0=0$.

Thus, $v$ satisfies an ultra-parabolic-like equation.
More precisely, 
it is not in general of the ultra-parabolic type since $\tilde a$ 
is not necessarily a positive definite matrix (although 
it satisfies the non-degeneracy assumption given in (A1)).

The main result of this section is given by the following theorem.
        
\begin{theorem}\label{T2}
Let $u$ be a quasi-solution to \eqref{d-p} under assumption (A1). 
Then there exists a subinterval $(\alpha',\beta')\subseteq (\alpha,\beta)$, 
$\alpha'<\beta'$, 
such that the function $s_{\alpha',\beta'}(u)$ admits the strong trace 
in the sense of Definition \ref{traces}. 
\end{theorem}

We shall prove this result by applying Lemma \ref{equivalence}. 
Here the (local) block structure \eqref{eq:atilde-block} of the 
diffusion matrix is important since it dictates the choice of 
the scaling \eqref{scal}, i.e.~the value of parameter $k$ 
will be precisely equal to $k$ given in (A1). 
The strong convergence of \eqref{scal} will be obtained using 
a compactness result of \cite{HKMP}, but to apply it 
the full non-degeneracy condition should hold.
Therefore, we first establish the result in the case when non-degeneracy 
conditions are satisfied also with respect to the first-order terms
(flux):
\begin{itemize}
\item[(A2)] For $\alpha,\beta$ given by (A1), there exists no interval 
$(\alpha_1,\beta_1)\subseteq(\alpha,\beta)$, $\alpha_1<\beta_1$, such that for some
$(\xi_0,\bar{\mxi})\in\R^{1+(d-k)}\setminus\{0\}$ the mapping 
$\lambda\mapsto \xi_0 + 
\langle \bar{f}(\lambda)\,|\,\bar{\mxi}\rangle$
is equal to zero on the whole interval $(\alpha_1,\beta_1)$, where
$f=\mff'$ and $f(\lambda)=(\tilde{f}(\lambda),\bar{f}(\lambda))\in
\R^k\times\R^{d-k}$.
\end{itemize}
Assumption (A2) will finally be removed by the inductive argument 
as it was done in \cite{pan_jhde}.
\begin{remark}
\label{R1}
We note that with (A1)--(A2) we require that neither of the mentioned 
mappings is not equal to zero on non-degenerate intervals, 
unlike the situation in Definition \ref{def:non-deg-ab} 
where the mappings should 
not be zero simultaneously on non-degenerate intervals. 
However, if the matrix $a$ has (locally) form \eqref{eq:atilde-block},
then these conditions are equivalent, since in (A1)--(A2) 
we have separated the directions of degeneracy of the flux and the diffusion.   
\end{remark}

\begin{theorem}\label{sc-ndeg}
Let $u$ be a quasi-solution to \eqref{d-p} under assumptions (A1)--(A2). 
Then the function $v=s_{\alpha,\beta}(u)$ admits the strong trace 
in the sense of Definition \ref{traces}.
\end{theorem}
\begin{proof}
For $k$ given in (A1), let us use the notation of the previous section:
$\hat{\mx}=(\tilde{\mx},\bar{\mx})\in\R^k\times\R^{d-k}$, $\my=(\tilde\my,\bar{\my})\in\R^k\times\R^{d-k}$ and 
$(u_m)$ given by \eqref{scal}.
By Lemma \ref{equivalence} and Remark \ref{rem:sab-sol} it is 
sufficient to prove that the sequence
$$
v_m(\hat{t},\hat{\mx},\my):=s_{\alpha,\beta}\bigl(u_m(\hat{t},\hat{\mx},\my)\bigr)
$$
strongly converges in $L^1_{loc}(\R^{d+1}_+\times\Rd)$.

By rescaling equation \eqref{eq:sab_e-c} we get that $v_m$ satisfies 
the same equation as the one obtained for $u_m$ in the previous 
section, except that terms \eqref{eq:unbdd_terms_non-diag} and 
\eqref{eq:unbdd_terms_diag} do not appear (since $v_m\in [\alpha,\beta]$ 
and $A$ satisfies \eqref{eq:atilde-block}), i.e.~for any
$\lambda\in\R$ and $\my\in\Rd$ it holds
\begin{align}
\nonumber
\pa_{\hat{t}} |v_m -\lambda|
	&+\Div_{\bar\mx}\Bigl({\rm sgn}(v_m-\lambda)
	\bigl(\bar\mff(v_m)-\bar\mff(\lambda)\bigr)\Bigr)
	-\sum\limits_{s,j=1}^k \pa_{\tilde{x}_j \tilde{x}_s}
	\Bigl({\rm sgn}(v_m-\lambda) 
	\bigl({A}_{sj}(v_m)-{A}_{sj}(\lambda)\bigr)\Bigr)\\
&= -\frac{1}{m}\tilde{\gamma}_m^\my(\hat{t},\hat{\mx},\lambda) 
	-\frac{1}{\sqrt{m}}\Div_{\tilde\mx}\Bigl({\rm sgn}(v_m-\lambda)
	\bigl(\tilde\mff(v_m)-\tilde\mff(\lambda)\bigr)\Bigr)\,=:\eta_m^\my \;.
\label{new}
\end{align}  
i.e. $v_m$ is also a quasi-solution to \eqref{d-p}.

Clearly, for every $\lambda \in \R$ and $\my\in \R^d$, the sequence $\frac{1}{\sqrt{m}}\Div_{\tilde\mx}\Bigl({\rm sgn}(v_m-\lambda)
	\bigl(\tilde\mff(v_m)-\tilde\mff(\lambda)\bigr)\Bigr)$ 
converges strongly to zero in $H^{-1}_{loc}(\R^{d+1}_+))$. 
Moreover, $\frac{1}{m}  \tilde{\gamma}_m^\my$ converges to zero in 
${\cal M}(\R^{d+1}_+)$ (this is proved in the same way as in
\cite[Lemma 3]{vass} or \cite[Lemma 3.2]{pan_jhdeB}). 
Since  ${\cal M}(\R^{d+1}_+)$ is compactly embedded in 
$W^{-1,r}_{loc}(\R^{d+1}_+)$, $r\in [1,\frac{d+1}{d})$, we 
finally see that for any $\lambda\in\R$ and $\my\in\Rd$ 
\begin{equation}\label{eq:compact-tmp}
\pa_{\hat{t}} |v_m -\lambda|
	+\Div_{\bar\mx}\Bigl({\rm sgn}(v_m-\lambda)
	\bigl(\bar\mff(v_m)-\bar\mff(\lambda)\bigr)\Bigr)
	-\sum\limits_{s,j=1}^k \pa_{\tilde{x}_j \tilde{x}_s}
	\Bigl({\rm sgn}(v_m-\lambda) 
	\bigl({A}_{sj}(v_m)-{A}_{sj}(\lambda)\bigr)\Bigr)
\end{equation}
is compactly embedded in $W^{-1,r}_{loc}(\R^{d+1}_+)$, $r\in [1,\frac{d+1}{d})$. 
Moreover, since for any $\alpha_1,\beta_1\in \R$,
$\alpha_1<\beta_1$, and any function $g:\R\to\R^l$ it holds
\begin{equation*}
g(s_{\alpha_1,\beta_1}(w))= 
	\frac{g(\alpha_1)+g(\beta_1)}{2}+\frac{{\rm sgn}(w-\alpha_1)
	\bigl(g(w)-g(\alpha_1)\bigr)-{\rm sgn}(w-\beta_1)
	\bigl(g(w)-g(\beta_1)\bigr)}{2} \;,
\end{equation*}
compactness of \eqref{eq:compact-tmp} implies that for 
any $\alpha_1,\beta_1\in\R$, $\alpha_1<\beta_1$, and any $\my\in\Rd$
\begin{equation}\label{ab-cond}
\pa_{\hat t}s_{\alpha_1,\beta_1}(v_m)+\Div_{\bar{\mx}}\bar\mff
	\bigl(s_{\alpha_1,\beta_1}(v_m)\bigr)
	-\sum\limits_{s,j=1}^k \pa_{\tilde{x}_j \tilde{x}_s}A_{sj}
	\bigl(s_{\alpha_1,\beta_1}(v_m)\bigr)
\end{equation} 
is compactly embedded in $W^{-1,r}_{loc}(\R^{d+1}_+)$, 
$r\in [1,\frac{d+1}{d})$, as well.
Indeed, one just needs to 
apply the identity above for $g(w)=w$, $g=\mff$ and $g=A$, 
and apply it to \eqref{eq:compact-tmp}.

Therefore, keeping in mind that (A1)--(A2) imply that the non-degeneracy
conditions on $(\alpha,\beta)$ are satisfied in the sense of Definition
\ref{def:non-deg-ab}, we can apply 
Lemma \ref{lm:compactness} below 
to conclude that for every $\my\in \R^d$ there exists 
a subsequence of $(v_m(\cdot,\my))$ (not relabelled) 
strongly converging in $L^1_{loc}(\R^{d+1}_+)$ 
toward say $\tilde{v}_0(\cdot,\my)$.

On the other hand, from \eqref{new} (see also Remark \ref{rem:sab-sol})
we have for ${h}^{\alpha,\beta}_m={\rm sgn}(v_m-\lambda)$:
\begin{align}
\label{kin-1}
 &\pa_{\hat t} {h}^{\alpha,\beta}_m +\Div_{\bar{\mx}}  
	\bigl( \bar{f}(\lambda) \, {h}^{\alpha,\beta}_m \bigr) -
	D^2_{\tilde{\mx}} \bigl(a(\lambda)  {h}^{\alpha,\beta}_m\bigr) = 	\pa_\lambda  \eta_m^\my \,, \\
\label{kin-2}
&{h}^{\alpha,\beta}_m{\big|}_{\hat t=0}={h}^{\alpha,\beta}_0 \,,
\end{align} where ${h}^{\alpha,\beta}_0$ is the weak trace provided 
by Theorem \ref{w-t} of the corresponding kinetic function
$h^{\alpha,\beta}={\rm sgn}(v-\lambda)$. 
If we denote $h^{\alpha,\beta}_m\overset{\scriptscriptstyle{m\to \infty}}{\longrightarrow} {\tilde h}^{\alpha,\beta}:={\rm sgn}(\tilde{v}_0(\cdot,\my)-\lambda)$  (the limit is subsequential), letting $m\to \infty$ in \eqref{kin-1}--\eqref{kin-2} along the chosen subsequence, we see that ${\tilde h}^{\alpha,\beta}$ satisfies the Cauchy problem
\begin{align}
\label{lim-1}
&\pa_{\hat t} {\tilde h}^{\alpha,\beta} +\Div_{\bar{\mx}} 
	\bigl( \bar{f}(\lambda) \,{\tilde h}^{\alpha,\beta} \bigr) -	D^2_{\tilde{\mx}}  (a(\lambda)  {\tilde h}^{\alpha,\beta}) = 0 \,,\\
\label{lim-ic}
&{\tilde h}^{\alpha,\beta}\Big|_{t=0}={h}^{\alpha,\beta}_0(\my,\lambda) \,,
\end{align} which implies ${\tilde h}^{\alpha,\beta}\equiv {h}^{\alpha,\beta}_0(\my,\lambda)$ since the solution of the latter Cauchy problem is unique. Moreover, from here it follows that the entire sequence $(\int_{\R} {h}^{\alpha,\beta}_m(\cdot,\lambda)\rho(\lambda) d\lambda)$ must converge (strongly) toward $(\int_{\R} {h}^{\alpha,\beta}_0(\my,\lambda)\rho(\lambda) d\lambda)$ (since every converging subsequence must converge toward the solution to \eqref{lim-1}--\eqref{lim-ic}).

Thus, we see that 
\begin{equation}
\label{ent-kin}
\int_{\R} {h}^{\alpha,\beta}_0(\my,\lambda)\rho(\lambda) d\lambda=\int_{\R}{\rm sgn}\bigl(\tilde{v}_0(\cdot,\my)-\lambda\bigr)\rho(\lambda) d\lambda \,,
\end{equation} for every $\rho\in C^1_c(\R)$ and any $\my\in\Rd$. Since $\tilde{v}_0$ is bounded between $-M$ and $M$, where $M:=\{|\alpha|,|\beta|\}$, we can take $\rho(\lambda)=\chi_{[-M,M]}(\lambda)$ -- the characteristic function of the interval $[-M,M]$. Inserting such $\rho$ in \eqref{ent-kin}, we get
$$
\tilde{v}_0(\hat{t},\hat{\mx},\my)=\frac{1}{2}\int_{-M}^M {h}^{\alpha,\beta}_0(\my,\lambda)\rho(\lambda) d\lambda \,,
$$
i.e.~$\tilde{v}_0$ does not depend on $\hat{t}$ and $\hat{x}$ and 
it is equal to the weak trace of $v=s_{\alpha,\beta}(u)$ (see \eqref{eq:u_w-t}).
The convergence $v_m(\cdot\,,\my)\xrightarrow{m\to\infty} \tilde{v}_0(\my)$ in $L^1_{loc}(\R^{d+1}_+)$ holds along entire sequence. According to the Lebesgue dominated convergence theorem, this actually means that $(v_m)$ converges in $L^1_{loc}(\R_{\hat{t},\hat{\mx}}^{d+1}\times \R^d_{\my})$ (i.e.~with respect to all three variables $\hat{t},\hat{\mx}$, and $\my$).

We have thus proved that $(v_m)=(s_{\alpha,\beta}(u_m))$ satisfies conditions of Lemma \ref{equivalence} and this in turn implies that any quasi-solution to \eqref{d-p} under assumptions (A1)--(A2) indeed admits existence of strong traces.
\end{proof}

In the following lemma we extract from  
\cite[Corollary 27]{HKMP}
a compactness result suitable for our application.

\begin{lemma}\label{lm:compactness}
Let $(u_m)$ be a bounded sequence in $L^\infty(\R^{d+1}_+)$ such that
$u_m(t,\mx)\in [\alpha,\beta]$, a.e. $(t,\mx)\in\R^{d+1}_+$. 
Let $r>1$ be such that for any $\alpha_1,\beta_1\in [\alpha,\beta]$,
$\alpha_1<\beta_1$, the sequence
\begin{equation*}
\pa_{t}s_{\alpha_1,\beta_1}(u_m)+\Div_{{\mx}}\mff
	\bigl(s_{\alpha_1,\beta_1}(u_m)\bigr)
	- D^2_\mx\cdot A\bigl(s_{\alpha_1,\beta_1}(u_m)\bigr)
\end{equation*} 
is precompact in $W^{-1,r}_{loc}(\R^{d+1}_+)$, where 
$\mff\in C^1(\R; \Rd)$ and $A\in C^1(\R;\R^{d\times d})$ satisfy the
non-degeneracy conditions on $(\alpha,\beta)$ in the sense of Definition \ref{def:non-deg-ab}. 

Then $(u_m)$ contains a subsequence convergent in $L^1_{loc}(\R^{d+1}_+)$. 
\end{lemma}

The only important difference between this lemma and \cite[Corollary 27]{HKMP}
is that in \cite[Definition 2]{HKMP} it is required that 
the non-degeneracy conditions hold on the whole $\R$. 
However, since $u_m$ takes values only in the segment $[\alpha,\beta]$,
it is enough to assume that the non-degeneracy conditions hold for every 
subinterval $(\alpha_1,\beta_1)\subseteq (\alpha,\beta)$. 
Indeed, in the notation of \cite[Theorem 25]{HKMP}, it is easy to 
see that $\mu^{pp}=0$ for $p\not\in [\alpha,\beta]$, where $\mu^{pp}$ 
is an H-measure given in \cite[Proposition 11]{HKMP} 
(cf.~\cite{Ger, Tar}). Thus, it is left only to 
show that the H-measure $\mu^{pp}$ is equal to zero also for $p\in[\alpha,\beta]$,
for which the non-degeneracy conditions are needed. 
Another way to see this is by noting that 
we can smoothly extend $\mff$ and $A$ out of the interval $(\alpha,\beta)$ 
so that the non-degeneracy conditions hold globally.

\smallskip

Let us now prove Theorem \ref{T2}.

\noindent \textbf{Proof of Theorem \ref{T2}:}
We shall use the method of induction with respect to the space dimension $d$ as introduced in \cite{pan_jhde} for $A=0$.
Let $\alpha,\beta,k$ be as in (A1).

In the case when $k=d$ (i.e.~when there is no flux part of the equation), 
we can use Theorem \ref{sc-ndeg} to infer about existence of strong traces for 
$s_{\alpha,\beta}(u)$.

Assume now that equation \eqref{d-p} is given on $d$-dimensional space and that in this situation any quasi-solution to \eqref{d-p} admits the strong trace. 
We shall prove from here that a quasi-solution $u$ satisfying \eqref{d-p} in $d+1$-dimensional case also satisfies the conclusion of Theorem \ref{T2}.

If condition (A2) is fulfilled, then we use Theorem \ref{sc-ndeg} to conclude
that the statement of Theorem \ref{T2} holds, since the inductive step is made. 
If (A2) fails to hold, then
there exist $(\xi_0,\bar\mxi)\in\R^{1+(d+1-k)}\setminus\{0\}$ and 
a non-degenerate interval $(\alpha_1,\beta_1)\subseteq (\alpha,\beta)$ such that
\begin{equation}
\label{flux-cond}
\xi_0+\sum\limits_{j=1}^{d+1-k} f_{k+j}(\lambda)\bar\xi_j = 0 \;, 
	\quad \lambda\in(\alpha_1,\beta_1) \,.
\end{equation}
It is obvious that it must be $\bar\mxi\neq 0$,
so the unit vector $\bar{\mxi}^*=\bar{\mxi}/|\bar\mxi|$ is well defined. 
Moreover, without loss of generality we can assume that
$\bar{\mxi}^*=(0,0,\dots,0,1)\in\R^{d+1-k}$,
since otherwise we just need to rotate the coordinate system
with respect to last $d+1-k$ coordinates, which will not affect
second order terms (as the first $k$ coordinates remain unchanged).
Thus, the following change of variables is regular
\begin{equation}\label{cov}
\bar{z}_{1}=\bar{x}_{1},\;\bar{z}_{2}=\bar{x}_{2},
	\dots, \;\bar{z}_{d-k}=\bar{x}_{d-k},
	\; \bar{z}_{d+1-k}= \xi_0 t+ \bar{x}_{d+1-k}
\end{equation} 
(we use $\mx=(\tilde{\mx},\bar{\mx})$). 
The equation \eqref{eq:sab-ultrapar-k-f} for $v=s_{\alpha_1,\beta_1}(u)$ 
in new variables 
becomes, for any $\lambda\in\R$, independent of the variation with respect 
to ${z}_d=\bar{z}_{d+1-k}$,
since
\begin{equation*}
\begin{aligned}
\pa_t \tilde h + \Div_{\tilde\mx} \bigl( \tilde f(\lambda)\, \tilde h\bigr) 
	+ \sum\limits_{j=1}^{d-k} \pa_{\bar{z}_j}\bigl(f_{k+j}(\lambda)
	\,\tilde h\bigr)
	&+\pa_{\bar{z}_{d+1-k}}\Bigl(\bigl(\xi_0+ f_{d+1}(\lambda)
	\bigr) \,\tilde{h}\Bigr) \\
&-D^2_{\tilde \mx}\cdot\bigl( \tilde a(\lambda) \,  \tilde h\bigr) 
	= \pa_\lambda
	\bigl(\tilde \gamma(t,\tilde{\mx},\bar{\mz},s_{\alpha_1,\beta_1}(\lambda))\bigr) \;,
\end{aligned}
\end{equation*}
and \eqref{flux-cond} holds,
 where $\tilde{h}(t,\mx,\lambda)=h^{\alpha_1,\beta_1}(t,{\mx}',t\xi_0+{x}_{d+1},\lambda)$, $\mx'=(x_1,x_{2},\dots,x_{d})$.
 Thus, for almost every $z_d=\bar{z}_{d+1-k}$ on $\R^d_+\times\R$
 it holds
 \begin{equation*}
 \pa_t \tilde h + \Div_{\tilde\mx} \bigl( \tilde f(\lambda)\, \tilde h\bigr) 
 	+ \sum\limits_{j=1}^{d-k} \pa_{\bar{z}_j}\bigl(f_{k+j}(\lambda)
 	\,\tilde h\bigr)
 -D^2_{\tilde \mx}\cdot\bigl( \tilde a(\lambda) \,  \tilde h\bigr) 
 	= \pa_\lambda
 	\bigl(\tilde \gamma(t,\tilde{\mx},\bar{\mz},s_{\alpha_1,\beta_1}(\lambda))\bigr) \;,
 \end{equation*}
  (see \cite[Section 6]{pan_jhde} for details in an analogical situation). 
  In this way, keeping in mind the equivalence between quasi- and kinetic solutions (see Theorem \ref{T7}), we have actually reduced the dimension of the equation and thus, we are in a position to use an induction argument with respect to dimension of the space. We note that, at this moment, we cannot claim existence of the strong traces on the entire interval $(\alpha,\beta)$ (see the first step in the induction argument) since the dimension of the equation is reduced only on the interval $(\alpha_1, \beta_1)$. This is however enough to conclude the statement of the theorem.
\endproof

\section{Existence of traces in the general case -- proof of Theorem \ref{main-traces}}\label{sec:main-traces}

In this section, we shall prove the main theorem of the paper -- Theorem \ref{main-traces}.
The proof is divided into four steps. 
\smallskip

\noindent \textbf{I.}
Let us take an arbitrary non-degenerate interval $(\alpha,\beta)\subseteq\R$.
We shall prove that there exists a subinterval on which we have existence of 
the strong trace, i.e.~that there exist $\alpha',\beta'\in\R$, $\alpha\leq \alpha'
<\beta'\leq \beta$, 
such that the function $s_{\alpha',\beta'}(u)$ admits the strong trace.

By Theorem \ref{T2} this is trivial in the case when condition (A1) 
is fulfilled on $(\alpha,\beta)$.

Thus, let us assume that condition (A1) does not hold on $(\alpha,\beta)$.
Then there exist a subinterval $(\alpha_1,\beta_1)\subseteq (\alpha,\beta)$, 
$\alpha_1<\beta_1$, and $\mxi^*\in \Sdmj$ such that for any 
$\lambda\in (\alpha_1,\beta_1)$ it holds
\begin{equation*}
\langle a(\lambda)\mxi^*\,|\,\mxi^* \rangle=0 \;.
\end{equation*} We can rotate the coordinate system if necessary in order to change the coordinates in the frame of which $\mxi^*=e_d=(0,\dots,0,1)$
(the flux in new coordinates would still be a $C^1$ functions, which 
is the only requirement that we need on this term). Assuming the latter, notice that $a(\lambda)=\sigma(\lambda)\sigma(\lambda)$ for the symmetric positive 
semi-definite matrix $\sigma$ (the principle square root), from where we have
$$
\langle a(\lambda)\mxi^*\,|\,\mxi^* \rangle=\langle \sigma(\lambda)\mxi^*\,|\, \sigma(\lambda)\mxi^* \rangle=0 \,, \ \ \lambda \in (\alpha_1,\beta_1) \,.
$$ 
Since $\mxi^*=e_d$, from the above we get that for any 
$j\in\{1,2,\dots,d\}$ it holds 
$$
\sigma_{jd}(\lambda)=0 \,, \ \ \lambda \in (\alpha_1,\beta_1) \,.
$$ 
Thus, for every $j=1,\dots,d$,
$$
a_{jd}(\lambda)=\sum_{k=1}^d \sigma_{jk}(\lambda)\sigma_{kd}(\lambda)=0 \,, 
	\ \  \lambda \in (\alpha_1,\beta_1) \,,
$$ 
i.e.~matrix $a$ in the interval $(\alpha_1,\beta_1)$ has the block form 
\eqref{eq:atilde-block} for $k=d-1$. 

If the whole assumption (A1) is satisfied on $(\alpha_1,\beta_1)$,
i.e.~if the matrix $\tilde{a}$ satisfies non-degeneracy assumption 
on $(\alpha_1,\beta_1)$, Theorem \ref{T2} is applicable.

If not, we continue with the procedure above. 
It is obvious that in a finite number of steps we reach to
a non-degenerate interval $(\alpha_l,\beta_l)\subseteq(\alpha,\beta)$
on which (A1) is fulfilled.
Indeed, maximal number of steps is equal to the dimension of the space,
as in that case we get (in notation of (A1)) $k=0$.
Thus, we can apply Theorem \ref{T2} on $(\alpha_l,\beta_l)$ 
to get the claim.

We have thus proved that 
\begin{equation}
\label{property}
\begin{aligned}
&\textit{for every interval $(\alpha,\beta)$ there exists a subinterval $(\alpha',\beta')$}\\
&\textit{such that $s_{\alpha',\beta'}(u)$ admits the strong trace.}
\end{aligned} 
\end{equation}

\noindent \textbf{II.} 
Let $M:=\|u\|_{L^\infty(\R^{d+1}_+)}$, i.e.~$u(t,\mx)\in [-M,M]$ 
for almost every $(t,\mx)\in\R^{d+1}_+$.
Applying \eqref{property}, we shall construct an open dense subset 
$I$ of $[-M,M]$ with a property that for any $(\alpha,\beta)\subseteq I$
function $s_{\alpha,\beta}(u)$ admits the strong trace.

Let us consider all (non-degenerate) open intervals with rational endpoints 
contained in $(-M, M)$ and arrange them in a sequence denoted by $(\tilde I_j)$.
Let us denote by $I_j\subseteq \tilde I_j$ an open interval given by
\eqref{property}.
We define $I=\bigcup_j I_j$. 

Set $I$ is clearly open and contained in $[-M,M]$. 
Moreover, it is dense in $[-M,M]$.
Indeed, if that were not the case, then there would be an open 
interval $I'\subseteq [-M,M]$ such that $I'\cap I=\emptyset$.
However, for some $j\in\N$ it holds $I_j\subseteq \tilde I_j\subseteq I'$,
thus $\emptyset\neq I_j\subseteq I'\cap I$ leads to a contradiction.

Furthermore, it is easy to see that the set obtained from $I$ by
removing countable number of points is dense in 
$[-M,M]$ as well.

By Remark \ref{rem:sab-sol},
for any $(\alpha,\beta)\subseteq I$ we have 
that $s_{\alpha,\beta}(u)$ admits the strong trace,
i.e.~we can consider any subinterval of $I$,
not necessarily equal to $I_j$. 

Let us note that although $I$ is open and dense, set $[-M,M]\setminus I$
could still be of strictly positive Lebesgue measure
(e.g.~$[-M,M]\setminus I$ might be a (fat) Cantor set of strictly positive 
Lebesgue measure).
\smallskip

\noindent \textbf{III.}
Let us take a sequence $(t_m)$ in $E$ (see \eqref{eq:set_E})
converging to zero and denote by 
$$u_m:=u(t_m,\cdot\,) \,.$$
By Remark \ref{rem:w-t_u} sequence $(u_m)$ converges 
weakly-$\star$ in $L^\infty(\Rd)$ to the weak trace $u_0$. 
In order to prove that $u$ admits the strong trace, by Remark 
\ref{rem:w-t_u_B} it is sufficient to prove that 
a subsequence of $(u_m)$ strongly converges in $L^1_{loc}(\Rd)$.  

Since 
\begin{equation}
\label{property-for-cnv}
\int_{-M}^M H(u_m-\lambda) \,d\lambda=u_m+M \,,
\end{equation} 
where $H$ is the Heaviside function,
it is enough to prove that $\bigl(H(u_m-\lambda)\bigr)_m$ converges 
strongly in $L^1_{loc}(\Rd)$ along a 
subsequence for almost every $\lambda\in [-M,M]$.

For a fixed $\lambda\in[-M,M]$, sequence 
of functions $H(u_m-\lambda)$ is bounded in $L^\infty(\Rd)$, 
hence converges weakly-$\star$ along 
a subsequence. A connection between this limit and the weak limit 
of $(u_m)$ can be expressed 
in terms of the corresponding \emph{Young measure} \cite{Young}, 
or more precisely using \emph{measure-valued functions}
\cite{Tar79, pan_matsb}, which we briefly present here in our setting
following aforementioned references.

A measaure-valued function on $\Rd$ is a weakly measurable mapping
$\Rd\ni\mx\mapsto\nu_\mx$ into the space of Borel probability measures
having compact supports on $\R$.
Weak measurability of $\nu_\mx$ means that for any continuous 
function $\ph$ on $\R$ the mapping 
$\mx\mapsto \int \ph(\lambda)\, d\nu_\mx(\lambda)$ is 
Lebesgue measurable on $\Rd$. 

There exist a subsequence of $(u_m)$ (not relabelled)
and a measure-valued function $\nu_\mx^0$ such that 
for any $\ph\in C(\R)$ we have
\begin{equation}\label{eq:mv-weak-conv}
\ph(u_m(\mx)) \xrightharpoonup[m\to\infty]{\;\;*\;\;} \int \ph(\lambda)\,d\nu_\mx^0(\lambda)
	\qquad \hbox{weakly-$\star$ in } L^\infty(\R^d_\mx) \;.
\end{equation}
Moreover, $\operatorname{supp}\nu_\mx^0\in [-M,M]$,
for a.e.~$\mx\in\Rd$. 

If we denote by $\delta_{\lambda_0}$ the Dirac measure at 
$\lambda_0\in\R$, then by
$\nu_\mx^m:=\delta_{u_m(\mx)}$
it is given a sequence of \emph{regular} measure-valued functions. 
Then, $\ph(u_m)$ can be expressed as 
\begin{equation*}
\ph(u_m) = \int \ph(\lambda) \,d\nu_\mx^m(\lambda)
\end{equation*} 
and \eqref{eq:mv-weak-conv} is often referred to as a \emph{weak} 
convergence of $(\nu_\mx^m)$ to $\nu_\mx^0$.  
It is not difficult to see that $(u_m)$ strongly converges to 
$u_0$ in $L^1_{loc}(\Rd)$ if and only if 
$\nu_\mx^0=\delta_{u_0(\mx)}$. 
Furthermore, if in \eqref{eq:mv-weak-conv} 
the weak-$\star$ convergence in $L^\infty(\Rd)$ is replaced by 
the strong convergence in $L^1_{loc}(\Rd)$, which is denoted 
by the \emph{strong} convergence of $(\nu_\mx^m)$ to $\nu_\mx^0$,
then $(u_m)$ strongly converges in $L^1_{loc}(\Rd)$ to $\int\lambda \,d\nu_\mx^0(\lambda)$.
In this case, by the uniqueness of the limit, it holds
$u_0(\mx)=\int\lambda \,d\nu_\mx^0(\lambda)$, i.e.~$\nu_\mx^0=\delta_{u_0(\mx)}$.

Although it will not be used in the following analysis, 
one might find interesting to notice that by Theorem \ref{w-t}
we have $\nu_\mx^0=-\frac{1}{2}\partial_\lambda h_0(\mx,\lambda)$
in distributional sense.

For $\mx\in\Rd$ and $\lambda\in \R$ we define
\begin{align*}
U_m(\mx,\lambda) &:= \nu_{\mx}^m\bigl((\lambda,+\infty)\bigr)
	= H\bigl(u_m(\mx)-\lambda\bigr) \,, \ m\in\N \,,\\
U_0(\mx,\lambda) &:= \nu_{\mx}^0\bigl((\lambda,+\infty)\bigr) \,.
\end{align*}
For any $m\in\N_0$ and $\lambda\in\R$ we have
$U_m(\cdot\,,\lambda)\in L^\infty(\Rd)$ (cf.~\cite{pan_matsb}), 
$0\leq U_m(\cdot\,,\lambda)\leq 1$, and 
for a fixed $\mx\in\Rd$ mappings 
$\lambda\mapsto U_m(\mx,\lambda)$, $m\in\N_0$, are monotonically
decreasing.
Moreover, since $M=\|u_m\|_{L^\infty(\Rd)}$ (see step II), 
it is easy to see that $U_m(\cdot\,,\lambda)=0$ for $\lambda>M$ and 
$U_m(\cdot\,,\lambda)=1$ for $\lambda<-M$, 
which is another evidence that only values of $\lambda$ that 
matters are within the set $[-M,M]$.

Let us define a set 
\begin{equation*}
P:= \Bigl\{ \lambda\in\R : \lim_{\lambda'\to\lambda}
	U_0(\cdot\,,\lambda')=U_0(\cdot\,,\lambda) \ \hbox{in} \
	L^1_{loc}(\Rd)\Bigr\} \;.
\end{equation*}
By \cite[Lemma 4]{pan_matsb} the complement $\R\setminus P$ is at most 
countable and for any $\lambda\in P$ we have
$U_m(\cdot\,,\lambda)\xrightharpoonup{\;*\;} U_0(\cdot\,,\lambda)$
weakly-$\star$ in $L^\infty(\Rd)$. 

Denote $V_m^\lambda (\mx) := U_m(\mx,\lambda)-U_0(\mx,\lambda)$. 
Thus, by the above, we have for any $\lambda\in P$
\begin{equation}\label{eq:Vlambda-conv}
V_m^\lambda \xrightharpoonup[m\to\infty]{\;\;*\;\;} 0
	\qquad \hbox{weakly-$\star$ in } L^\infty(\Rd) \;.
\end{equation}
Under this new notation, it is left to prove that (along a 
subsequence) $(V_m^\lambda)_m$ converges strongly in 
$L^1_{loc}(\Rd)$ to zero for almost every $\lambda\in[-M,M]$.
Indeed, since any continuous function can be approximated uniformly 
on any compact set by finite linear combinations of functions of the form 
$\lambda\mapsto H(\lambda-\lambda_0)$, this would imply that 
(along a subsequence) $(\nu_\mx^m)$ converges strongly to $\nu_\mx^0$.

Let $\lambda, \lambda'\in P$ and take a compact set $K\subseteq \Rd$. 
Since $U_m$ is monotone in variable $\lambda$, we have
\begin{align*}
\int_{K} \bigl|U_m(\mx,\lambda) &- U_m(\mx,\lambda')\bigr| \,d\mx
	= \biggl|\int_{K} 
	\bigl(U_m(\mx,\lambda) - U_m(\mx,\lambda')\bigr) \,d\mx\biggr| \\
&=\biggl|\int_{K} V_m^\lambda(\mx) \,d\mx
- \int_{K} V_m^{\lambda'}(\mx) \,d\mx
+ \int_{K} \bigl(U_0(\mx,\lambda)-U_0(\mx,\lambda') \bigr) 
	\,d\mx\biggr| \,.
\end{align*}
Using this identity we get
\begin{align*}
\int_{K} |V_m^\lambda(\mx)&-V_m^{\lambda'}(\mx)| \,d\mx \\
&\leq \int_{K} \bigl|U_m(\mx,\lambda) - U_m(\mx,\lambda')\bigr| \,d\mx
	+ \int_{K} \bigl|U_0(\mx,\lambda)-U_0(\mx,\lambda')\bigr|\,d\mx \\
&\leq \biggl|\int_{K} V_m^\lambda(\mx) \,d\mx\biggr|
+\biggl|\int_{K} V_m^{\lambda'}(\mx) \,d\mx\biggr|
+ 2\int_K\bigl|U_0(\mx,\lambda)-U_0(\mx,\lambda')\bigr|\,d\mx \;.
\end{align*}
Thus, since $\lambda,\lambda'\in P$ by \eqref{eq:Vlambda-conv} 
we obtain
\begin{equation*}
\limsup_{m\to\infty} \int_{K} |V_m^\lambda(\mx)-V_m^{\lambda'}(\mx)| 
	\,d\mx
\leq 2\int_K\bigl|U_0(\mx,\lambda)-U_0(\mx,\lambda')\bigr|\,d\mx \;.
\end{equation*}
Finally, since $\lambda\in P$ we have
\begin{equation}\label{eq:Vlambdalambdaprime}
\lim_{P\ni\lambda'\to \lambda} \limsup_{m\to\infty} 
\|V_m^\lambda-V_m^{\lambda'}\|_{L^1(K)} = 0 \;.
\end{equation}

In the next and final step of the proof we shall prove the strong convergence 
of $(V_m^\lambda)$ for almost every $\lambda$, i.e.~for 
$\lambda\in [-M,M]\cap P$.
\smallskip

\noindent \textbf{IV.}
We shall consider first the case $\lambda\in I\cap P$. 
This means that there exist an interval $(\alpha,\beta)\subseteq I$
such that $\lambda\in (\alpha,\beta)$ and the function 
$s_{\alpha,\beta}(u)$ admits the strong trace (see steps I and II).
In particular, the sequence $\bigl(s_{\alpha,\beta}(u_m)\bigr)$
converges strongly in $L^1_{loc}(\Rd)$
and let us denote the limit by $u_0^{\alpha,\beta}$.

Since it holds (recall that $\lambda\in (\alpha,\beta)$)
\begin{align*}
U_m(\cdot\,,\lambda) = H(u_m-\lambda) 
	= H(s_{\alpha,\beta}(u_m)-\lambda) \,,
\end{align*}
by passing to a subseqeunce (not relabelled) we have
\begin{equation*}
\lim_{m\to\infty} U_m(\cdot\,,\lambda) = H(u_0^{\alpha,\beta}-\lambda)
\end{equation*}
in $L^1_{loc}(\Rd)$.
Since the weak and the strong limits must coincide, we have 
$H(u_0^{\alpha,\beta}-\lambda) = U_0(\cdot\,,\lambda)$, 
i.e.~$V_m^\lambda$ converges strongly to zero in 
$L^1_{loc}(\Rd)$.

Let us take now $\lambda\in [-M,M]\cap P$. 
Since $I\cap P$ is dense in $[-M,M]$ (see step III), 
there exists a sequence $(\lambda_n)$ in $I\cap P$ such that 
$\lim_n\lambda_n=\lambda$.
Then for an arbitrary compact $K\subseteq \Rd$ we have
$$
\|V_m^\lambda\|_{L^1(K)} \leq \|V_m^\lambda-V_m^{\lambda_n}\|_{L^1(K)}
	+ \|V_m^{\lambda_n}\|_{L^1(K)} \,.
$$
Letting first $m$ to infinity and then $n$ to infinity, 
by \eqref{eq:Vlambdalambdaprime} and since 
$\bigl(V_m^{\lambda_n}\bigr)_m$ converges strongly to zero, 
we get that $\bigl(V_m^\lambda\bigr)_m$ converges strongly to zero as well. 

Thus, for almost every $\lambda\in[-M,M]$ we have
$$
\lim_{m\to\infty} U_m(\cdot\,,\lambda) = U_0(\cdot\,,\lambda)
$$
strongly in $L^1_{loc}(\Rd)$. 
Therefore, 
\begin{equation*}
u_m = -M + \int_{-M}^M H(u_m-\lambda)d\lambda
\end{equation*}
converges strongly in $L^1_{loc}(\R^{d+1}_+)$.
Thus, the proof is complete.
\endproof

\smallskip

\begin{remark}\label{rem:weak-strong-trace}
In Step IV od the proof of Theorem \ref{main-traces} we denoted 
by $u_0^{\alpha,\beta}$ the strong trace of $s_{\alpha,\beta}(u)$.
However, at the end we can easily conclude that it must be 
$u_0^{\alpha,\beta}=s_{\alpha,\beta}(u_0)$, since \emph{the whole}
function $u$ admits the strong trace (see Remark \ref{rem:sab-sol}).

In fact, this can be concluded even if we know only that 
$s_{\alpha,\beta}(u)$ has the strong trace.
Indeed, since $(s_{\alpha,\beta}(u_m))$ converges strongly in 
$L^1_{loc}(\Rd)$ to $u_0^{\alpha,\beta}$ and 
the function $\lambda\mapsto s_{\alpha,\beta}(\lambda)$ is continuous, 
by \eqref{eq:mv-weak-conv}
it holds that for a.e.~$\mx\in\Rd$ and any $\varphi\in C(\R)$ we have
$$
\int \varphi(s_{\alpha,\beta}(\lambda))\,d\nu_\mx^0(\lambda)
	= \int \varphi(\lambda) \,d\delta_{u_0^{\alpha,\beta}(\mx)}(\lambda)
	= \varphi(u_0^{\alpha,\beta}(\mx)) \,.
$$
By the arbitrariness of $\varphi$ (e.g.~we can choose it to be equal to zero 
on $\langle \alpha,u_0^{\alpha,\beta}(\mx)\rangle$
or $\langle u_0^{\alpha,\beta}(\mx),\beta\rangle$)
and continuity of $s_{\alpha,\beta}$,
this implies that 
$s_{\alpha,\beta}(\lambda)=u_0^{\alpha,\beta}(\mx)$
on $\operatorname{supp}\nu_\mx^0$, for a.e.~$\mx\in\R^d$
(a more general result in this direction can be found
in \cite[Corollary 2.6]{pan_jhdeB}).
Since $s_{\alpha,\beta}$ is monotone, the 
last identity can be extended to the (closed) convex hull of 
$\operatorname{supp}\nu_\mx^0$. 
It is left to notice that the \emph{barycentre} 
$\int \lambda\,d\nu(\lambda)$ of any 
probability measure $\nu$ is contained in the closed convex hull of its support,
hence
$$
u_0^{\alpha,\beta}(\mx) 
	= s_{\alpha,\beta}\left(\int\lambda\,d\nu_\mx^0(\lambda)\right)
	= s_{\alpha,\beta}(u_0) \,,
$$
for a.e.~$\mx\in\Rd$. 

We would like to thank to the referee for this remark. 
\end{remark}

Let us close this section with a proof of Corollary \ref{main-traces-cor}.

\noindent \textbf{Proof of Corollary \ref{main-traces-cor}:}
Let $p>1$ and let $u\in L^\infty_{loc}(\R^+;L^p_{loc}(\R^{d}))$ be a quasi-solution 
to \eqref{d-p}. 
By Theorem \ref{main-traces} (see Remark \ref{rem:sab-sol})
for any $n\in\N$ the function $s_{-n,n}(u)$ admits the strong trace
denoted by $u_0^n$. 
Thus, we can define a unique (up to equality a.e.) measurable function 
$u_0:\Rd\to [-\infty,+\infty]$ such that 
$s_{-n,n}(u_0)=u_0^{n}$ a.e.~in $\Rd$
(see Remark \ref{rem:sab-sol}).
It is left to prove that $u_0\in L^1_{loc}(\Rd)$. 

Let us take a sequence $(t_m)$ in $E$ (see \eqref{eq:set_E})
converging to zero and denote $u_m:=u(t_m,\cdot\,)$.
By the assumption $(u_m)$ is a bounded sequence in $L^p_{loc}(\Rd)$.
Then, for any $n\in\N$ we have that $s_{-n,n}(u_m)$
strongly converges to $u_0^n=s_{-n,n}(u_0)$ in $L^1_{loc}(\Rd)$
(here we used that the set \eqref{eq:set_E} does not depend on $n$;
see Remark \ref{rem:sab-sol}).
From here, using e.g.~\cite[Theorem 28]{HKMP}, we see that $(u_m)$ strongly converges to $u_0$ in $L^1_{loc}(\Rd)$,
which ensures $u_0\in L^1_{loc}(\Rd)$. 
\endproof

\section{Concluding remarks}\label{sec:conclusion}

To conclude the paper, we shall comment on possible extensions of the proved result.

\begin{itemize}

\item[1.] Following \cite{pan_jhde}, existence of traces can be proved under assumptions of mere continuity of the flux $\mff$ and the diffusion matrix $A$. However, the proof in this case is technically more demanding and, since its essence is substantially the same as the one presented here, we omit it.

\item[2.] Existence of strong traces remains to hold if the flux $\mff$ depends on $(t,\mx)$ (so called heterogeneous or non-autonomous flux), 
in which case in Definition \ref{def:quasisol} relation 
\eqref{e-c} modifies as follows:
\begin{align*}
\pa_t |u-\lambda| &+\Div_\mx\Bigl({\rm sgn}(u-\lambda)\bigl(\mff(t,\mx,u)-\mff(t,\mx,\lambda)\bigr)\Bigr)\\
&-D^2_\mx\cdot\Bigl( {\rm sgn}(u-\lambda) \bigl(A(u)-A(\lambda)\bigr)\bigr) 
= -\gamma(t,\mx,\lambda) \,.
\end{align*}
However, we need to assume the non-degeneracy conditions:
  
\emph{for a.e.~$\mx\in \R^{d}$ and for all $(\xi_0,\mxi) \in \Sd$
there is no interval $(\alpha,\beta)\subseteq\R$ such that
the functions  
\begin{equation*}
F(\lambda)= \xi_0+\langle \mff'(0,\mx,\lambda)\,|\,\mxi \rangle \ \ {\rm and} \ \ G(\lambda)=\langle A'(\lambda) \mxi\,|\,\mxi \rangle, \ \ \lambda \in (\alpha, \beta) \,,
\end{equation*} are identically equal to zero simultaneously on 
$(\alpha,\beta)$.}

The proof in this case would be analogous to the proof of Theorem \ref{sc-ndeg}. We have just a slight modification when introducing the blow-up change of variables \eqref{blow-up}.
Indeed, in \eqref{new} (for the non-autonomous flux $\mff$) 
one should rewrite
$$
\bar\mff\Bigl(\frac{\hat{t}}{m},\frac{\tilde{\mx}}{\sqrt{m}}
+\tilde{\my},\frac{\bar{\mx}}{m}+\bar{\my}, u_m\Bigr)
- \bar\mff\Bigl(\frac{\hat{t}}{m},\frac{\tilde{\mx}}{\sqrt{m}}
+\tilde{\my},\frac{\bar{\mx}}{m}+\bar{\my}, \lambda\Bigr)
$$ 
as
\begin{equation}\label{eq:non-auto-case}
\begin{aligned}
\Bigl(\bar{\mff}(0,\my,u_m)-\bar{\mff}(0,\my,\lambda)\Bigr)
&+ \biggl(\bar\mff\Bigl(\frac{\hat{t}}{m},\frac{\tilde{\mx}}{\sqrt{m}}
+\tilde{\my},\frac{\bar{\mx}}{m}+\bar{\my}, u_m\Bigr)
- \bar{\mff}(0,\my,u_m)\biggr) \\
&- \biggl(\bar\mff\Bigl(\frac{\hat{t}}{m},\frac{\tilde{\mx}}{\sqrt{m}}
+\tilde{\my},\frac{\bar{\mx}}{m}+\bar{\my}, \lambda\Bigr)
- \bar{\mff}(0,\my,\lambda)\biggr) \,.
\end{aligned}
\end{equation}

The terms corresponding to the last two summands we
put on the right-hand side of \eqref{new}, which read
\begin{equation*}
\begin{split}
\Div_{\bar\mx} & \left( 
{\rm sgn}(u_m-\lambda)
\left(
\bar\mff\Bigl(\frac{\hat{t}}{m},\frac{\tilde{\mx}}{\sqrt{m}}+\tilde{\my},\frac{\bar{\mx}}{m}+\bar{\my}, u_m\Bigr)
-\bar\mff(0,\tilde{\my},\bar{\my}, u_m) 
\right)\right)\\
&-\Div_{\bar\mx}\left({\rm sgn}(u_m-\lambda)\left(\bar\mff\Bigl(\frac{\hat{t}}{m},\frac{\tilde{\mx}}{\sqrt{m}}+\tilde{\my},\frac{\bar{\mx}}{m}+\bar{\my}, \lambda\Bigr)-\bar\mff(0,\tilde{\my},\bar{\my}, \lambda) \right) \right) \,,
\end{split}
\end{equation*}
but these terms converge to zero strongly in $H^{-1}_{loc}(\R^{d+1}_+)$ for every fixed $\hat{\my}=(\bar{\my},\tilde{\my})$.
Hence, these terms do not affect the proof procedure.

On the other hand, the term corresponding to the first summand 
in \eqref{eq:non-auto-case} remains on the left hand side 
(this is the reason why it is sufficient to have the non-degeneracy condition 
only for $t=0$).

\item[3.] If the diffusion depends on $(t,\mx)$, i.e.~if $A=A(t,\mx,\lambda)$, 
then in same special situations we can still get the result without 
non-degeneracy assumptions. 
More precisely, one needs that the vector $\mxi^*$ from step I of the proof 
of Theorem \ref{main-traces} is valid for a.e.~$(t,\mx)\in\R^d_+$.
For example, for the diffusion matrix $a=A'$ of the form
$$
a(t,\mx,\lambda)=\begin{bmatrix}
\tilde{a}(t,\mx,\lambda)& 0\\
0 & \bar{a}(\lambda)
\end{bmatrix} \,,
$$
where $\tilde{a}(t,\mx,\lambda) \in \R^{k\times k}$ is a positive definite matrix,
the statement of Theorem \ref{main-traces} still holds. 
The first $k$-coordinates of the flux $\mff$ in this case can also depend on $(t,\mx)$. If entire flux depends on $(t,\mx)$, then we need to assume the non-degeneracy conditions as in item 2.

In this case, the proof needs an additional argument of regularity of $u$ in the sense that $u_{x_j}\in L^2(\R^{d+1}_+)$, $j=1,2,\dots,k$ (see \cite{AM_jhde}).  
\end{itemize}

\noindent {\bf Acknowledgements.} 
The authors would like to thank two anonymous referees for
	their insightful comments that helped to improve the presentation 
	and the quality of the paper.

This work was supported in part by the Croatian Science Foundation under 
projects UIP-2017-05-7249 (MANDphy) and IP-2018-01-2449 (MiTPDE), and by the projects 
P33594 and M 2669 Meitner-Programm of the Austrian Science Fund FWF.

Permanent address of D.~Mitrovi\'c is University of Montenegro.

\end{document}